\documentclass[12pt]{article}
\usepackage[utf8]{inputenc}
\usepackage{amsmath, amssymb, amsthm}

\usepackage{comment}
\usepackage{tabularx}
\usepackage{enumerate,float,subcaption,multirow}
\usepackage[dvipsnames]{xcolor}
\usepackage{hyperref,tikz, tikzpagenodes}
\hypersetup{colorlinks = true, citecolor = black, linkcolor = NavyBlue, urlcolor = NavyBlue}
\usetikzlibrary{through,intersections}
\usepackage{algpseudocode,algorithm}
\usepackage[sorting=nyt,sortcites=true,maxbibnames=9,autopunct=true,autolang=hyphen,hyperref=true,abbreviate=false, backend=biber,giveninits=true,doi=false,url=false]{biblatex}
\renewbibmacro{in:}{} 
\DeclareFieldFormat{journaltitle}{{#1}} 
\DeclareFieldFormat{booktitle}{{#1}} 
\DeclareFieldFormat[article,incollection,inproceedings,thesis]{title}{\mkbibemph{#1}} 
\DeclareFieldFormat{pages}{#1} 
\usepackage{mathrsfs}
\usepackage{tikz-cd}
\usepackage{cleveref}

\colorlet{ss}{black!20}
\colorlet{sslabel}{black!50}
\colorlet{backline}{black!10}

\newtheorem{theorem}{Theorem}
\newtheorem{corollary}[theorem]{Corollary}
\newtheorem{lemma}[theorem]{Lemma}

\newtheorem*{remark}{Remark}

\usepackage[left=30mm,right=30mm,top=20mm,bottom=20mm]{geometry}

\def\0{{\bf 0}}
\def\1{{\bf 1}}

\def\In{{\mathrm{In}}}

\def\Q{\mathbb{Q}}

\newcommand{\F}{{\mathbb {F}}}
\newcommand{\Z}{{\mathbb{Z}}}

\newcommand{\tr}{{\,\mathrm{tr}\,}}

\newcommand{\erase}[2][says]{}
\makeatletter
\newcommand{\rref}[2]{\hyperref[#2]{{#1}~\ref*{#2}}}
\makeatother

\title{The limits of inertia bounds on graph powers}
\author{Aida Abiad\thanks{\texttt{a.abiad.monge@tue.nl}, Department of Mathematics and Computer Science, Eindhoven University of Technology, The Netherlands; Department of Mathematics and Data Science, Vrije Universiteit Brussel, Belgium} \and Nils van de Berg\thanks{\texttt{n.p.v.d.berg@tue.nl}, Department of Mathematics and Computer Science, Eindhoven University of Technology, The Netherlands} \and Sam Mattheus\thanks{\texttt{ sam.mattheus@vub.be}, Department of Mathematics and Data Science, Vrije Universiteit Brussel, Belgium}}

\date{}

\begin{document}

\maketitle

\begin{abstract}

We investigate the tightness of inertia and inertia-type bounds related to graph powers. This class of graphs appears in several contexts, such as the theory of distance-regular graphs, quantum information theory and coding theory. They are a priori not covered by previous separation constructions involving locally sparse graphs, see Kwan and Wigderson (2024), or Alon's pseudorandom graphs, see Elphick, Tang and Zhang (2025). Instead, we apply the latter's strategy using complements of pseudorandom graphs, improved by Chen and Zhang (2026), to random Cayley graph constructions, whose graph-power structure is ensured by additive combinatorial arguments. From this we are able to conclude that also for graph powers there can be a large gap between inertia bounds and the actual independence number. For the inertia-type bound, which is based on polynomials, we achieve even larger separation by considering sequences of graphs whose spectrum converges to a continuous distribution. Our methods draw from several areas of mathematics, ranging from additive combinatorics and probabilistic methods, to character sums over finite fields and arithmetic geometry.\\

\end{abstract}

\section{Introduction}

The independence number of a graph $G$, denoted by $\alpha(G)$, is one of the most studied parameters in graph theory. Since it is NP-hard to compute \cite{Karp72} and to approximate \cite{Zuckerman06}, a lot of effort has been put into finding good bounds on the independence number. The spectral approach has yielded many such bounds. The \textit{spectrum} $\mathrm{Spec}(G)$ of a graph $G$ is the multiset of the eigenvalues of its adjacency matrix, which can be ordered as $\lambda_1 \ge \lambda_2 \ge \ldots \ge \lambda_n$. 
Two of the most well-known types of spectral bounds are ratio bounds (that use the ratio between the largest and smallest eigenvalues), and inertia bounds (that count the signs of the eigenvalues). While the ratio bounds are relatively well understood \cite{haemershoffman, filmus2022EKR, zhou_unified_2025}, the inertia bounds are much less well understood.

The first inertia bound was introduced in its unweighted form by Cvetkovi\'c \cite{cvetkovic_graphs_1971} in 1971. This bound was strong enough to give an algebraic proof of the Erd\H{o}s-Ko-Rado theorem, which can also be obtained using the ratio bound, see e.g.\ \cite{godsilmeagherbook}. However, the full power of the bound lies in the freedom to choose weights. Calderbank and Frankl \cite{Calderbank_Frankl_1992} studied the weighted version of the inertia bound in 1992 in order to derive multiple results in extremal set theory. For a Hermitian matrix $M$, we let $n_{\geq 0}(M)$ be the number of non-negative eigenvalues of $M$ (counted with multiplicity). For a graph $G$ a \emph{weighted adjacency matrix} of $G$ is a Hermitian matrix $B$ with zero diagonal where $B_{ij} = 0$ when there is no edge $i \sim j$ in $G$.

\begin{theorem}[Weighted inertia bound, cf. \cite{Calderbank_Frankl_1992}]\label{thm:weightedinertia}
For a graph $G$ and a weighted adjacency matrix $B$, 
    \[
    \alpha(G)  \leq n_{\geq 0}(B).
    \]
\end{theorem}

 We let $n_{\geq 0}(G)$ denote the optimal bound on $\alpha(G)$ obtained using Theorem~\ref{thm:weightedinertia}. The strength of this bound is such that it is tight for every graph on at most $10$ vertices \cite{Elzinga2010}. The search space of weights is infinite and there is no practical program known to find optimal weights. There are some heuristics known, such as for Cayley graphs in \cite{Rooney2014cocliques}, but no method that guarantees optimality. This also makes it hard to lower bound $n_{\ge 0}(G)$ by anything larger than $\alpha(G)$. 

In 2018 Sinkovi\'c published the first example where using real weights cannot give a matrix such that $n_{\ge 0}(B)$ is equal to the independence number \cite{sinkovic2018}. Only in the last two years has the separation between $\alpha(G)$ and $n_{\ge 0}(G)$ been stretched to arbitrarily large gaps. Kwan and Wigderson used ideas from probability theory to show that $C_4$-free graphs have a weighted inertia bound that grows linearly in the order $n$ \cite{kwan2024}. Such graphs can have ratio bound and hence independence number bounded by $O(n^{3/4})$. Elphick, Tang and Zhang \cite{tang2025} showed that there is no way to upper bound $n_{\ge 0}(G)$ in terms of $\alpha(G)$ by demonstrating a family of graphs with fixed independence number, $\alpha(G)=2$, where the weighted inertia bound grows at a rate $\Omega(n^{1/9})$. In this paper we continue their line of work and show similar separation results for the more general $k$-independence number of a graph. 

The \emph{$k$-independence number of a graph} $G$, denoted $\alpha_k(G)$, is a natural generalization of the independence number that considers sets where all vertices are at distance greater than $k$ from each other.  This more general independence parameter is also NP-hard to compute, see  \cite{Kong1993OnSets}. Equivalently, it can be defined as the independence number of $G^k$, the \emph{$k$-th graph power} of $G$, the graph on the same vertex set as $G$ where two vertices are adjacent if they are at distance at most $k$ in $G$. 

Graph powers help design efficient algorithms for specific combinatorial optimization problems, see \cite{AB1995,BC2000}. In distributed computing, the $k$-th power of graph $G$ represents the possible flow of information during $k$ rounds of communication in a distributed network of processors organized according to a graph $G$, see \cite{l92}. The investigation of the behavior of the independence number of graph powers is motivated by questions in (quantum) information theory and coding theory. Indeed, $\alpha_k(G)$ and hence inertia bounds can be used to estimate the possible size of minimum distance $k+1$ codes, see e.g.\ \cite{APR2025, reijnders2025}. The inertia can also be used to upper bound the quantum $k$-independence number \cite{WEA2022} and the quantum distance-$k$ chromatic number \cite{AM2025,AJ2025}, which is particularly helpful because these parameters are not known to be computable. 

Good inertia bounds typically require a good choice of weights. For graphs of interest this is much more difficult to find for $G^k$, than for $G$, because the graph is denser and thus has more freedom in its choice of weights. This led Abiad, Coutinho and Fiol \cite{kindep2019} to introduce a more restricted weight optimization version of the inertia bound, which bounds $\alpha_k(G)$ directly in terms of the spectrum of $G$. Note that the relation between $G$ and $G^k$ is not particularly algebraic, so their spectra are in general not related, see e.g.\ \cite{abiad_optimization_2022,Das2013LaplacianGraph,JZ2017}. Therefore bounds on $\alpha_k$ which use the spectrum of $G$ directly (such as the following bound) carry different information than those using the spectrum of $G^k$. 

\begin{theorem}[Inertia-type bound \cite{kindep2019}] \label{thm:inertia-type}
    Let $P_k$ denote the vector space of real polynomials in $x$ with degree at most $k$. Let $p \in P_k$ and $\lambda_1 \geq \lambda_2\geq \cdots \geq \lambda_n$ be the eigenvalues of the adjacency matrix $A$ of $G$. Let $w(p) = \min_i (p(A))_{ii}$ and $W(p) = \max_i (p(A))_{ii}$ denote the minimum and maximum of the diagonal of $p(A)$ and $$\In_k(p,G) := \min\{|\{i : p(\lambda_i) \geq w(p)\}|, |\{i : p(\lambda_i) \leq W(p)\}|\}.$$
    Then $\alpha_k(G) \leq \In_k(p,G)$.
\end{theorem}

We write $\In_k(G)$ for the optimal inertia-type bound, i.e., using the optimal polynomial. Even less is known about the performance of this inertia-type bound. Finding $\In_k(G)$ can be formulated as a mixed-integer linear
program \cite{abiad_optimization_2022}, which can be solved in polynomial time for fixed $k$ \cite{Hoesel26}. Nevertheless, relatively little is known about the graphs for which the optimum equals $\alpha_k(G)$. Indeed, we do not know when the inertia-type bound is tight, except for one class of graphs attaining equality for $p = x^k$ \cite{ACT2016} and for specific Odd graphs \cite{abiad_optimization_2022}. 

In this paper we show that large separation exists between the $k$-independence number and inertia bounds. We show this for both the weighted inertia bound on $G^k$ and for the inertia-type bound on $G$. During the preparation of this paper, it was shown that the weighted inertia bound is at least as strong as the inertia-type bound \cite{notekindep2026}. Indeed, the separation we show for the inertia-type bound is larger. In Section~\ref{sec:weightedinertia} we use the approach from \cite{tang2025} to deduce separation results for the weighted inertia bound. Just after the appearance of the preprint version of this paper, Chen and Zhang managed to improve the results from \cite{tang2025} and we will use their stronger theorem \cite[theorem~4.1]{ChenZhang26}. With this we prove that certain non-explicit Cayley graphs have $\alpha_k(G) = 2$ and $n_{\geq 0}(G^k) = \tilde{\Omega}(n^{1/(2k)})$.

 \newtheorem*{thm:main1}{Theorem \ref{thm:kindepseparation}}
\begin{thm:main1}
    For fixed $k \geq 3$ there is a sequence of graphs $G_n$ of order $N$ such that $\alpha_k(G_n) = 2$ and $n_{\geq 0}(G_n^k) \geq \widetilde{\Omega}(N^{1/(2k)})$.
\end{thm:main1}

Separation also holds for $k = 2$, but our proof does not give the same $\frac{1}{2k}$ exponent. We cannot extend Kwan and Wigderson's result \cite{kwan2024} to get a linearly growing weighted inertia bound, because graph powers are not $C_4$-free and are generally locally dense. Instead, the proof of our separation result relies on the observation that all large enough subsets of certain Abelian groups can be written as a sumset $A + A$. This follows from an extension of Alon's result on this topic in \cite{alonsumsets2007} to symmetric sets and larger sums $kA$, see Theorem~\ref{thm:alonextend}, which may be a result of independent interest. A month before the appearance of this preprint, Chae and Lee posted an article that also obtained the extension of Alon's result to larger $k$ but without symmetry \cite{ChaeLee26}. In Section~\ref{sec:inertiatype} we show that the inertia-type bound can grow linearly with the order while the diameter is at most $k$ (implying $\alpha_k(G) = 1$). We show this for random graphs of the right density (Corollary~\ref{cor:inertiatypefar_sparse}) and inspired by the recent construction of pseudorandom graphs by Forey, Fres\'an, Kowalski and Wigderson in \cite{FFKW25} we also give explicit vertex-transitive graphs for which this holds.

 \newtheorem*{thm:main2}{Theorem \ref{thm:cayleyinertiatypefar}}
 \begin{thm:main2}
   Fix an integer $ k \ge 5$ and prime $p > k$. Define $s = \lfloor \frac{k}2\rfloor +1$. There exists an integer $a \ge 1$ such that $s \mid p^{2a} -1$. For any large enough prime power $q = (p^{2a})^h$ the graph 
    \[
    G_q^{s} := \mathrm{Cay}(\F_q^{s}, S_{s}), \quad \text{ where } S_{s} = \{\pm(t, t^2, \ldots, t^s) : t \in \F_q^\times\}
    \]
    has
    \[\alpha_k(G_q^{s}) = 1, \quad \quad \text{and} \quad \quad \In_k(G_q^{s}) = \Omega(|V(G_q^{s})|) \text{ as } q \to \infty.\]
  
\end{thm:main2}

The same separation result also holds for $2 \leq k \leq 4$, but requires different graphs and a different proof, which we only remark upon.

\subsubsection*{Notation}
The following two sections can be read independently. Therefore each of them will contain its own necessary preliminaries. In both sections we prove results for the Erd\H{o}s-R\'enyi random graph. We let $G(n,p)$ denote this random graph on $n$ vertices where each edge is included independently with probability $p$. The results we prove will be asymptotic. We use the following notation. For two positive functions $f,g : \mathbb{N} \to \mathbb{R}_{> 0}$ we write $f = o(g), \omega(g), O(g)$ or $\Omega(g)$ if the fraction $\frac{f(n)}{g(n)}$ as $n\to \infty$ goes to zero or infinity, is bounded from above or bounded from below by a constant $C > 0$ respectively. If $f = O(g)$ and $f = \Omega(g)$ we write $f = \Theta(g)$. We also use $f = \tilde{O}(g)$ if $f(n) = O(gh)$ where $h$ is a polylogarithmic function. The notation $\tilde{\Omega}$ is used analogously. We use $\log$ to denote logarithms with base 2 and $\ln$ for natural logarithms. Lastly, we say that a sequence of events $E_n$ occurs with high probability if $\mathbb{P}(E_n) = 1 - o(1)$ as $n \to \infty$.

\section{Separation of the weighted inertia bound}
\label{sec:weightedinertia}

Recall that $n_{\ge 0}(G)$ denotes the optimal weighted inertia bound for a graph $G$. In \cite{tang2025} the authors lower bounded  $n_{\geq 0}(G)$ for complements of pseudorandom graphs. This lower bound was recently improved by a factor $3$ in the exponent in \cite{ChenZhang26}. We use this improved lower bound to find separation results for the weighted inertia bound of graph powers. 
Let $\lambda(G)$ denote $\max\{|\lambda_2(G)|, |\lambda_n(G)|\}$, the second largest absolute value of an eigenvalue. 

\begin{theorem}[\cite{ChenZhang26}]\label{thm:expanderinertia}
    Suppose $G = (V,E)$ is an $n$-vertex graph, such that its complement contains a spanning $d$-regular subgraph $\Gamma$ with $0 < d < n-1$ and eigenvalues $d = \lambda_1 \ge \lambda_2 \ldots \ge \lambda_n$. Then $\lambda_2 > \lambda_n$ and $$n_{\geq 0}(G) \geq \frac{d -\lambda_n}{\lambda_2 -\lambda_n}.$$
    In particular, if $\lambda > 0$ satisfies $|\lambda_i| \le \lambda$ for every $i \ge 2$, then
    \[
    n_{\ge 0}(G) \ge \frac{d+ \lambda}{2\lambda}.
    \]
\end{theorem}

In \cite{kwan2024} Kwan and Wigderson conjectured that the weighted inertia bound is typically $\Omega\left(\frac{n}{\log n}\right)$. A lower bound of $\Omega(\left(n/\log n\right)^{1/2})$ follows from Theorem~\ref{thm:expanderinertia} and the fact that a random graph is pseudorandom. The full conjecture got resolved through different means in \cite{ChenZhang26}.


 Theorem~\ref{thm:expanderinertia} was used to show that the complements of Alon's construction of triangle-free pseudorandom graphs \cite{alonRamsey1994} are examples of graphs with $\alpha(G) = 2$ and $n_{\geq 0}(G) = \Omega(n^{1/3})$. In Section~\ref{sec:inertiagraphpowers} we show a similar separation result for graph powers. 

It can be shown that the graphs arising as the $k$th power of a graph form only a small fraction of all graphs and hence it seems reasonable that the separation should become smaller as $k$ increases. We could not prove that the complements of Alon's triangle-free graphs from \cite{alonRamsey1994} are graph powers for any $k \geq 2$ and large $n$. Therefore we could not rely on these graphs to give separation results. Instead, we rely on random Cayley graphs to give graphs $G^k$ with large separation. These graphs still have unbounded weighted inertia bound while $\alpha_k(G) = 2$. To show that they arise as graph powers we first need to extend Alon's result on sumsets, specifically Theorem~2.1 in \cite{alonsumsets2007}, which we extend in Section~\ref{sec:ksumsets}.

\subsection{Large sets in finite vector spaces are k-sumsets}\label{sec:ksumsets}

Define $kA$ to be the set of sums of exactly $k$ elements of $A$, allowing repetitions, i.e. the set $\{a_1 + \ldots + a_k : a_i \in A\}$. In this section we extend Alon's sufficient condition for a set to be a sumset \cite{alonsumsets2007} to $k > 2$. In particular, we prove the following:

\begin{theorem}\label{thm:alonextend}
    Fix an integer $k \geq 2$. Let $q > k$ be a prime number, $\Gamma = (\Z/q\Z)^n$, and $N := |\Gamma|$. For $N$ sufficiently large, there exists a constant $D = D(k) > 0$ such that for any subset $F \subseteq \Gamma$ of size at most $D\sqrt[k]{\frac{N}{(\log N)^{k-1}}}$ there is a set $A \subseteq \Gamma$ such that $kA = \Gamma \setminus F$.
\end{theorem}

The same result with an explicit constant  was recently proven in \cite[Theorem~1.3]{ChaeLee26}. They stated their result for $F$ of order $N^{1/k}/\log N$, but their proof also supplies $N^{1/k}/\log^{(k-1)/k}N$. In order to apply Theorem~\ref{thm:alonextend} to random Cayley graphs we need extra conditions on $A$. Firstly, we require that $A$ is symmetric, so that $G = \mathrm{Cay}(\Gamma, A\setminus0)$ is an undirected graph and $A - A = A+A$. Secondly, we require that $0 \in A$ so that $G^k = \mathrm{Cay}(\Gamma, kA \setminus 0 )$. Hence, we will prove a statement with slightly more conditions. A similar, but simpler, proof (without the need for orbits, Lemma~\ref{lem:noBfullorbit} and exceptions for $0$) would also prove Theorem~\ref{thm:alonextend}. 

\begin{theorem}\label{thm:alonextendsym0}
    Fix an integer $k \geq 2$. 
    Let $q > k$ be a prime number, $\Gamma = (\Z/q\Z)^n$, and $N := |\Gamma|$. For $N$ sufficiently large, there exists a constant $D = D(k) > 0$ such that for any symmetric subset $F \subseteq \Gamma$ with $0 \notin F$ and of size at most $D\sqrt[k]{\frac{N}{\log^{k-1} N}}$ there is a symmetric set $A \subseteq \Gamma$ with $0 \in A$ such that $kA = \Gamma \setminus F$.
\end{theorem}

Throughout the proof of Theorem~\ref{thm:alonextendsym0} we will assume that $N$ is large with respect to the fixed $k$. All asymptotic notation is henceforth with regard to $N \to \infty$, where the implicit constants in big O notation will depend on $k$. 

The proof largely follows the same structure as the proof by Alon in \cite{alonsumsets2007}. We will show that conditioned on four likely events, the natural generalization of the random deletion process works. However, extra technicalities arise because of the symmetry, the forced inclusion of $0$ and because for $k > 2$ there are distinct sets that are not disjoint, but whose sums are equal. We first give the full proof and then show that each of the ``bad" events $E_0, E_1^g,E_2^g,E_3^g$ happen with negligible probability in the lemmas that follow.

\begin{proof}
   Given $F$ and $k$ we let $B$ be the set of \emph{bad} multisets, i.e., those multisets of $\Gamma$ of size $k$ that sum to an element in $F$. For any (multi)set $V \subseteq \Gamma$ we write $\mathcal{O}(V)$ to be the (multi)set of orbits under the inversion action, i.e. $\{\{x,-x\} : x \in V\}$. We create a random subset $C \subseteq \Gamma$ by adding orbits $\mathcal{O}(x) = \{x,-x\}$ independently with probability 
   \begin{align}\label{eq:defineprob}
       p = (a_k\cdot2^{k-1}k!N^{1-k}\ln N )^{1/k}
   \end{align}
   for some large enough constant $a_k$ we determine later and adding $0$ deterministically.
   Define $B_C$ as the subset of $B$ containing those multisets that are supported in $C$ and let $\bigcup B_C$ denote all $0 \neq x \in \Gamma$ that appear in some element of $B_C$. We define the first bad event $E_0$. \begin{description}
       \item[] $(E_0)$ There is an element of $B_C$ that contains both elements of a non-zero orbit.
   \end{description}
   In Lemma~\ref{lem:noBfullorbit} we show that $E_0$ holds with probability $o(1)$. Let $\bar{E}_0$ denote the complement event. 
   Let $A = C \setminus \bigcup B_C$.  We assumed that $0 \notin F$, so that the $k$-multiset $\{0,0,\ldots,0\} \notin B$. Therefore every multiset of $B$ contains some non-zero element of $\bigcup B_C$, so no multiset of $B$ is supported by $A$. 
   Additionally, $F$ is symmetric, so $b \in B_C$ implies $-b \in B_C$ and hence $\bigcup B_C$ is symmetric. Thus $A$ is symmetric and we have defined it to contain $0$, so $kA$ will also contain $0$. We will show that each non-zero element in $\Gamma \setminus F$ is in $kA$ with probability $1 - o(1/N)$. From this it follows by the union bound that  $kA = \Gamma\setminus F$ with  probability $1-o(1)$ and hence there is at least one such $A$.
   
   We say that a (multi)set of orbits sums to $h$ if we can choose an element of each orbit in the set such that the sum of those elements is $h$.

   Fix a non-zero element $g \in \Gamma \setminus F$. We define the bad event $E_1^g$ as follows:
   \begin{description}
       \item[] $(E_1^g)$ there is no set of disjoint $k$-sets of $C$ that sum to $g$ of size $2\ln N$.
   \end{description} Suppose $E_1^g$ fails and let $C_g$ be a family of $2\ln N$ disjoint $k$-sets of orbits summing to $g$. If $g \notin kA$ each set in $C_g$ needs to intersect a set in $B_C$. Consider pairs $(b,c) \in B_C \times C_g$ such that $\mathcal{O}(b) \cap c \neq \emptyset$. We will say that this is a \textit{blocking pair} or $b$ \textit{blocks} $c$. In writing down all blocking pairs each $b$ can appear at most $k$ times, because $b$ contains at most $k$ different elements and all $c \in C_g$ have disjoint orbits. Thus we can greedily pick a sequence $S$ of $\frac{2}{k}\ln N$ blocking pairs $(b_i,c_i)$ such that $\mathcal{O}(b_i) \cap c_j = \emptyset$ if $i < j$ (at each step we delete all $c_j$ that are blocked by the chosen $b_i$). Now consider the next bad event $E_2^g$:
   \begin{description}
       \item[] $(E_2^g)$ there is $\{-x, x\} \in \Gamma$ (with $x = 0$ possible), a sequence of $b_1,\ldots, b_{k+2} \in B_C$ and distinct $c_1,\ldots,c_{k+2} \in C_g$ such that all $b_i$ contain $x$ or $-x$ but none of the $c_i$ contain $\{-x,x\}$, $b_i$ blocks $c_i$ and $\mathcal{O}(b_i) \cap c_j = \emptyset$ for $i < j$.
   \end{description}
   
   Suppose that $E_2^g$ fails. For any $(b_i,c_i) \in S$ and $x \in b_i$ there are at most $k+1$ other $(b_j,c_j) \in S$ with $x \in b_j$ or in one of the orbits of $c_j$, because including $(b_i,c_i)$ there are at most $k+2$: at most one where an orbit of $c_j$ contains $x$ and at most $k+1$ more due to the failure of $E_2^g$. 
   By definition of $S$ any $(b_i,c_i)$ can intersect at most $(k+1)k$ $(b_j,c_j) \in S$ with $j < i$. If $g \notin kA$, then there must be $\frac{2}{k(k(k+1)+1)}\ln N$ pairs of $(b_i,c_i) \in S$ that are completely disjoint from each other. 
   This implies the third bad event $E_3^g$:
   \begin{description}
       \item[]  $(E_3^g)$ there are $\frac{2}{k^2(k+2)}\ln N$ disjoint blocking pairs $(b_i,c_i) \in B_C \times C_g$.
   \end{description}
   We have deduced that if $E_0, E_1^g$, $E_2^g$ and $E_3^g$ all fail, then $g \in kA$. In particular, by the union bound \[\mathbb{P}(kA \neq \Gamma\setminus F) \le \mathbb{P}(E_0) + \sum_{0 \neq g \in \Gamma \setminus F} \left[\mathbb{P}(E_1^g) + \mathbb{P}(E_2^g \cap \bar{E}_0) + \mathbb{P}(E_3^g \cap \bar{E}_0)\right]
   \]
    For every non-zero $g \in \Gamma \setminus F$ each of the events $E_1^g, E_2^g \cap \bar{E}_0, E_3^g \cap \bar{E}_0$ occurs with probability $o(1/N)$ by Lemmas~\ref{lem:many_sets},  \ref{lem:fewintersectingbodies} and \ref{lem:fewdisjointbodies} respectively.
   Hence the sum is $o(1)$ and thus $\mathbb{P}(kA = \Gamma\setminus F) = 1 - o(1)$. 
\end{proof}

We now prove that the events hold with appropriately small probability in the order $E_1^g, E_0, E_2^g, E_3^g$ (for $0 \neq g \in \Gamma\setminus F$).

\begin{lemma}\label{lem:many_sets}
    Let $E_1^g$ be the event that there are fewer than $2\ln N$ pairwise disjoint sets of $k$ orbits in $C$ that sum to $g$. Then $\mathbb{P}(E_1^g) = o(1/N)$.
\end{lemma}

\begin{proof}
    Note that $g \neq 0$. Let $\mathcal{S}_g$ be the family of $k$-sets of orbits (excluding the orbit $\{0\}$) that sum to $g$. Recall that this means that we can pick one element from each orbit so that the elements sum to $g$. Let $C_g'$ be the subset of $\mathcal{S}_g$ which only contains sets that are fully inside $C$. We will show that the set $C_g'$ contains a family of disjoint sets of size $s = 2\ln N$ with probability $1 - o(1/N)$.  

    Let $s = 2\ln N$. Let $\mathcal{M}$ denote a maximal family of disjoint sets in $C_g'$ and let $M = |\mathcal{M}|$. We will show that $\mathbb{P}(M < s) = o(1/N)$. Let $U = \bigcup_{c \in \mathcal{M}} c$ be the set of all orbits that are in one of the elements of $\mathcal{M}$. If $M < s$, then there are  $kM < ks$ orbits in $U$. If there is a set in $C_g'$ disjoint from $U$, then we can add it to the family, contradicting the maximality of $\mathcal{M}$. For any subset $S \subseteq \Gamma$ let $X_S$ denote the number of sets in $C_g'$ that are disjoint from $S$. The set $U$ satisfies the events that $X_U = 0$ and that $U \subseteq \mathcal{O}(C)$. These events are independent, as they depend on disjoint elements to be included in $C$. Therefore we find that
     $$\mathbb{P}(M < s) \le \sum_{i = 0}^{s-1} \sum_{\substack{U \subseteq \mathcal{O}(\Gamma),  |U| = ki,\\ \text{union of disjoint orbits in } \mathcal{S}_g}} \mathbb{P}(X_U = 0)\mathbb{P}(U \subseteq \mathcal{O}(C)).$$ 
     
     It is immediate to see that $\mathbb{P}(U \subseteq \mathcal{O}(C)) = p^{ki}$. In estimating $\mathbb{P}(X_U = 0)$ we first count the number of sets in $\mathcal{S}_g$ that are disjoint from $U$, denoted by $\mathcal{S}_g^U$. Given $k-1$ orbits of $S \in \mathcal{S}_g $ that are outside $U$ there is at least one choice of the final element (and hence final orbit) to make the set of orbits sum to $g$. This overcounts the total amount by a factor $k$, once for each final orbit added. There are $\frac{N-1}{2}$ orbits to choose from in total, so this suggests that there are at least $\frac{1}k\binom{(N-1)/2 - |U|}{k-1}$ such sets. If the last orbit is not zero, outside $U$ and distinct from the first $k-1$ orbits, then this indeed gives a set in $\mathcal{S}_g^U$. The number of cases where the last element fails to satisfy those conditions is $O(N^{k-2}(|U|+k))$, which is $o(1)$ of all cases because $|U| < 2k\ln N$. Hence the number of sets in $\mathcal{S}_g^U$ is at least $\frac{N^{k-1}}{2^{k-1}k!}(1-o(1))$. Hence the expected size of $C_g^U = \mathcal{S}_g^U \cap {C\choose k}$ is $$\mu_U \ge \frac{N^{k-1}}{2^{k-1}k!}(1-o(1))p^k = (a_k-o(1))\ln N.$$ 
     
     To apply Janson's inequality, specifically $\mathbb{P}(X_U = 0) \le \exp(-\mu_U + \Delta/2)$, we estimate 
     $$\Delta = \sum_{\substack{c, c' \in \mathcal{S}_g^U \text{distinct} \\ c \cap c' \neq \emptyset}} \mathbb{P}(c, c' \in C_g^U).$$

    The sets of orbits $c$ and $c'$ sum to $g$ and are distinct. They can intersect in at most $k-1$ elements. Suppose $c,c'$ intersect in $i$ orbits. There are at most $2^{k-1}\binom{(N-1)/2}{k-1}$ sets in $\mathcal{S}_g^U$, because fixing $k-1$ orbits gives at most $2^{k-1}$ final elements which give sum $g$ for some choice of representative per orbit. Given $c$ there are $\binom{k}{i}$ ways to choose orbits of $c'$ in the intersections. Choosing $k-1-i$ other elements leaves at most $2^{k-1}$ choices for the last element. Hence we can choose at most $2^{k-1}\binom{k}{i}\binom{(N-1)/2-k}{k-i-1}$ sets $c' \in C_g'$ that intersect $c$ in $i$ orbits for all $1 \leq i \le k-1$. Hence we bound 
    \begin{align*}\Delta \le \sum_{c \in \mathcal{S}_g^U} \sum_{i = 1}^{k-1} 2^{k-1}\binom{k}{i}\binom{(N-1)/2}{k-1-i}p^{2k-i} &\le  2^{k-1}\binom{(N-1)/2}{k-1}\sum_{i = 1}^{k-1} 2^{k-1} \binom{k}{i}\binom{(N-1)/2}{k-1-i}p^{2k-i}\\ &= \sum_{i = 1}^{k-1} O(N^{2k-2-i}(N^{1-k}\ln N)^{2-i/k}) = o(\ln N),
    \end{align*}
  
    so by Janson's inequality we have 
    $$\mathbb{P}(X_U = 0) \le \exp(-\mu_U + \Delta/2) \le \exp(-a_k\ln N +o(\ln N)) \leq N^{-a_k+o(1)}.$$ 
  
    The number of sets $U$ of size $ki$ that can arise as a union of disjoint sets in $\mathcal{S}_g$ is at most ${|\mathcal{S}_g| \choose i} \leq \left(\frac{e|\mathcal{S}_g|}{i}\right)^i \le \left(\frac{eN^{k-1}/k!}{i}\right)^i$, because $|\mathcal{S}_g| \le \frac{N^{k-1}}{k!}(1-o(1))$. Hence
    \begin{align*}
    \mathbb{P}(M < s) &\le \sum_{i = 0}^{s-1} \left(\frac{eN^{k-1}/k!}{i}\right)^ip^{ki}N^{-a_k+o(1)}\\
    &\le s\left(\frac{eN^{k-1}p^k}{s\cdot k!}\right)^sN^{-a_k+o(1)} \\
    &= s(e2^{k-2}a_k)^{2\ln N}N^{-a_k+o(1)} \\
    &= o(1/N),
    \end{align*}
    where $a_k$ is chosen large enough that this final equality holds. Here we used that $\left(\frac{eN^{k-1}p^k}{i\cdot k!}\right)^i$ is increasing in $i$ when $i \le \frac{N^{k-1}p^k}{k!}$ which is bigger than $2\ln N$ by definition of $p$, see \eqref{eq:defineprob}. 

    We conclude that the set $C_g'$ contains a family of disjoint sets of size $s = 2\ln N$ with probability $1- o(1/N)$. 
\end{proof}

Before we move to the next lemma, we make the following observation. Consider a $k$-multiset $b$ with all distinct elements fixed except one, which we give the variable $x$ and $x$ occurs $t$ times in $b$. We know that the sum of the elements of $b$ is $tx+h$ for some $h \in \Gamma$. We know that $t \leq k$, so $\gcd(t,q) = 1$ and thus there are exactly $|F|$ choices for $x$ such that $tx+h \in F$.

\begin{lemma}\label{lem:noBfullorbit}
    Let $E_0$ be the event that there is a multiset $b \in B_C$ that contains $\{x,-x\}$ for some non-zero $x$. Then $\mathbb{P}(E_0) = o(1)$.
\end{lemma}

\begin{proof}
    Each full orbit sums to zero, so in case $k = 2$ there are no sets containing $\{x,-x\}$ in $B$. Assume $k \ge 3$. If a multiset $b \in B_C$ contains a full non-zero orbit, then there must be a submultiset $b' \subseteq b$ without full orbits (also without $\{0\}$) of size at most $k-2 \ge 1$ that sums to an element in $F$. Let $X_t$ be the number of multisubsets of size $t$ of elements of orbits in $C$ containing no full orbits (including zero) and summing to an element in $F$. Suppose $b'$ of size $t$ contains no full orbits and contains $i$ different elements of $\Gamma$. Then $\mathbb{P}(\mathcal{O}(b') \subseteq C) = p^i$, because none of the elements lie in the same orbit and none are zero. There are at most $O(1)\binom{N}{i-1}|F|$ such $b'$, because there are at most $t^i$ ways to decide the multiplicities of the $i$ distinct elements and given those multiplicities choosing the first $i-1$ elements leave at most $|F|$ options for the final one, by the remark above. Therefore, the expected number of $X_t$ is
    \[
    \mathbb{E}(X_t) = \sum_{i = 1}^t O(1)N^{i-1}|F|p^i= \tilde{O}\left(\sum_{i = 1}^t N^{(i+1)/k-1}\right) = \tilde{O}_k(N^{(t+1-k)/k}).
    \]
    As $t \leq k-2$ all these expected values go to zero and the probability that there is one such set also goes to zero by Markov's inequality.
\end{proof}

In the remaining Lemmas we will count the number of blocking pairs given that the elements of $b$ come from a fixed set. We record this count here as a separate lemma.

\begin{lemma}\label{lem:blockingpairs}
    Consider the subset $B' \subseteq B$ of multisets that contain no non-zero $\{-x,x\}$ and let $S \subseteq \Gamma$ be a fixed set of size $s$. The number of blocking pairs $(b,c) \in B' \times \mathcal{S}_g$ with the elements of $b$ either in $S$ or one of the orbits of $c$ is  $O(|F|\binom{N}{k-2})$ where the implied constant only depends on $s$ and $k$.
\end{lemma}

\begin{proof}
    Let $z_1,\ldots,z_k$ be variables that denote the elements of the orbits of $c$ that sum to $g$. Without loss of generality let $z_1$ be the variable such that $\{z_1,-z_1\} \in \mathcal{O}(b)$, which exists because $b$ blocks $c$. If $z_1 \in b$, then, because $g \notin F$, there is a $z_i$ that will not appear in $b$. Similarly if $-z_1 \in b$, then there is a $-z_i$ that does not appear in $b$, because $-g \notin F$. Without loss of generality take $i = 2$. There are at most $\binom{N}{k-2}$ ways to choose $z_3,\ldots, z_k$ (less because they all need to be in distinct orbits). Leaving $z_1, z_2$ as variables there are at most $(s+2k)^k$ ways to decide on the values of $b\setminus \{z_1,-z_1,z_2,-z_2\}$  and the multiplicities of the $\pm z_1, \pm z_2$ that appear. Note that at most one of $z_1,-z_1$ appears in $b$, because $b \subseteq B'$. Let the multiplicity of $z_1$ be $t_1$ (where $t_1 < 0$ if $-z_1 \in b$). Similarly at most one of $z_2, -z_2$ appears and let $t_2$ be the (signed) multiplicity of $z_2$. Let $a = z_3+\ldots +z_k$ and let $a'$ be the sum of the remaining elements in $b$, then the conditions on $c$ and $b$ imply that
    \[
    \begin{pmatrix} 1 & 1\\ t_1 & t_2\end{pmatrix}\begin{pmatrix} z_1 \\z_2\end{pmatrix} = \begin{pmatrix} g - a\\ f - a' \end{pmatrix}
    \]
    where $f \in F$. We know that $|t_1| + |t_2| \le k < q$ and that $t_1,t_2$ do not have the same sign (note that $t_2$ can be zero). Hence the matrix has rank $2$ and for each $f \in F$ there is exactly one solution. We conclude that there are $O(|F|\binom{N}{k-2})$ blocking pairs.
\end{proof}

\begin{lemma}\label{lem:fewintersectingbodies}
     Let $E_2^g$ be the event that there is some orbit $\{x,-x\} \subseteq \Gamma$ and $m = k+2$ multisets $b_1,\ldots, b_m \in B_C$ and sets of orbits $c_1,\ldots,c_m \in C_g$  such that for all $i$ we have $\{x,-x\} \in \mathcal{O}(b_i)\setminus c_i$, $b_i$ blocks $c_i$, and $\mathcal{O}(b_i) \cap c_j = \emptyset$ for $i < j$. Then 
    $\mathbb{P}(E_2^g \setminus E_0) = o(1/N)$.
\end{lemma}

\begin{proof}
    We estimate the expected number of tuples without assuming that the $c_i$ are within $C_g$, i.e., satisfying the slightly looser conditions:  $x,b_1,\ldots, b_m, c_1,\ldots, c_m$ where $E_0$ does not hold, $b_i \in B_C$ blocks $c_i \in \binom{C}{k}$, $\{x,-x\} \in \mathcal{O}(b_i)\setminus c_i$ for all $i$, the $c_i$ are disjoint and sum to $g$, and $\mathcal{O}(b_i) \cap c_j = \emptyset$ if $i<j$. 
    
    First, we upper bound the number of such tuples with $b_i \in B$ and $c_i \in \mathcal{S}_g$. Suppose that the total number of elements in the set $L = \bigcup_{i = 1}^m\mathcal{O}(b_i)\setminus \bigcup_{i = 1}^m c_i$ is $\ell$. Then $\ell \le m(k-2)+1$ because each $\mathcal{O}(b_i)$ contains $\{x,-x\}$ and at least one element of $c_i$. There are at most $\binom{N}{\ell}$ choices for the $\ell$ orbits and $\mathcal{O}(x)$ must be chosen among those, so there are $\ell$ choices for $\mathcal{O}(x)$. In determining the number of possible tuples we fix the elements of $(b_i,c_i)$ in increasing order of $i$. Having determined $(b_j,c_j)$ for all $j < i$ we know that the elements of $b_i$ are among the elements in $L \cup \bigcup_{j < i} c_j$ a set of size less than $2(\ell +ki) = O(1)$, so by Lemma~\ref{lem:blockingpairs} we get at most $O(|F|\binom{N}{k-2})$ options for $(b_i,c_i)$. We can conclude that the number of tuples we consider is at most $O\left(\binom{N}{\ell}\ell\prod_{i = 1}^m \binom{N}{k-2}|F|\right)$. Among these there are $O\left(\binom{N}{\ell-1}\ell\prod_{i = 1}^m \binom{N}{k-2}(\ell+ik)^{k-2}|F|\right)$ sets with $0 \in L$.  
    In total there are at least $mk+\ell$ disjoint orbits in $\bigcup_{i = 1}^m\mathcal{O}(b_i) \cup \bigcup_{i = 1}^m c_i$, because all the $c_i$ are pairwise disjoint. The probability that a tuple is included is $p^{mk+\ell}$ or $p^{mk+\ell-1}$ in case $\{0\} \in L$. The case with $\{0\} \notin L$ has a factor $Np = \omega(1)$ more and is therefore the dominating term.
    Hence the expected number of such tuples is at most
    \begin{align*}
    O\left(N^\ell|F|^m\binom{N}{k-2}^mp^{mk+\ell}\right) &=  \tilde{O}(N^{\ell+m/k + m(k-2)-m(k-1)-\ell+\ell/k})\\ &= \tilde{O}_k(N^{(m+\ell)/k-m}) = o(1/N)
    \end{align*}
    
   \noindent because $m = k+2$ and $\ell \le m(k-2)+1$. By Markov's inequality we conclude that $\mathbb{P}(E_2^g\setminus E_0) = o(1/N)$.
\end{proof}

\begin{lemma}\label{lem:fewdisjointbodies}
      Let $E_3^g$ be the event that there are $\frac{2}{k^2(k+2)}\ln N$ disjoint blocking pairs $(b_i,c_i) \in B_C \times C_g$. Then,
    $\mathbb{P}(E_3^g \setminus E_0) = o(1/N)$.
\end{lemma}

\begin{proof}
    We only need to calculate the probability when none of the $b_i$ contain both elements of a non-zero orbit. For each $k \leq \ell \leq 2k-2$ define $A_\ell$ to be the event that there are at least $k+1$ disjoint blocking pairs $(b,c)$ such that as a set $|\mathcal{O}(b) \cup c| = \ell$ for every pair.
    Define $A_{2k-1}$ to be the similar event that there are at least $s := \frac{1}{k^2(k+2)}\ln N$ disjoint blocking pairs $(b,c)$ with $|\mathcal{O}(b)\cup c| = 2k-1$ for every pair.
    In order to show that $\mathbb{P}(E_3^g) = o(1/N)$, it suffices to show that $\mathbb{P}(A_\ell) = o(1/N)$ for all $k \leq \ell \leq 2k-1$ as long as $k(k-1)+s \leq \frac{2}{k^2(k+2)}\ln N$, which is clearly satisfied for large enough $N$.
    
   So consider for each $k \leq \ell \leq 2k-1$ the pairs $(b,c) \in B \times \mathcal{S}_g$ for which $|\mathcal{O}(b) \cup c| = \ell$. There are at most $\binom{(N+1)/2}{\ell -k}$ ways to choose the $\ell - k$ orbits that are not in $c$. 
   Given those, $b$ must take its elements among at most $2(\ell -k)$ elements in these orbits and the elements in $c$, so by Lemma~\ref{lem:blockingpairs} there are $O(\binom{N}{\ell-k}\binom{N}{k-2}|F|)$ blocking pairs with $|\mathcal{O}(b) \cup c| = \ell$. Hence there are at most $\binom{O(N^{\ell-2}|F|)}{t}$ choices for $t$ tuples $(b_i,c_i)$ with $|\mathcal{O}(b_i) \cup c_i| = \ell$ and $\mathcal{O}(b_i) \cap c_i \neq \emptyset$. If we require that the $(b_i,c_i)$ are pairwise orbit-disjoint, then the probability that all their elements are in $C$ is $p^{t\ell}$ unless zero is in one of the $b_i$. Without loss of generality $b_1$ has the zero, which gives a factor $\Omega(N)$ less in the number of possible $b_1$, a factor $t$ extra in the binomial and a factor $p$ less in the probability that $b_1 \subseteq C$. This gives a factor $O(t/(Np))$ extra. This factor goes to zero, so this case has even smaller probability. Therefore there is a constant $C_k$ only dependent on $k$ such that the probability that $A_\ell$ holds is at most

    {\footnotesize{
    \begin{align}\binom{O(N^{\ell-2}|F|)}{t}p^{t\ell}  &\le \left(C_k\frac{N^{\ell-2}|F|p^\ell}{t}\right)^t \notag\\ 
    &= \left(C_k\frac{N^{-2+\ell/k}(\ln N)^{\ell/k}|F|}{t}\right)^t. \label{eq:Et3}
    \end{align} 
    }}
    For $\ell \leq 2k-2$, we have $t = k+1$ 
     and $|F| = \widetilde{O}(N^{1/k})$, so this is $\tilde{O}(N^{-\frac{1}{k}})^{k+1} = o(\frac1N)$. For $\ell = 2k-1$, $t = s = \frac{1}{k^2(k+2)}\ln N$, and $|F| \le DN^{1/k}(\ln N)^{-(k-1)/k}$ we get that

    \begin{align*}\mathbb{P}(A_{2k-1}) &\leq \left(C_k \cdot D \cdot k^2(k+2)\right)^s = o\left(\frac1N\right),
    \end{align*}
    for sufficiently small $D$.
\end{proof}

\begin{remark}
    The condition $k(k-1) + s \le 2s$ implies that the result requires that $N = 2^{\Omega(k^5)}$. This can be slightly improved to $N = k^{\Omega(k^4)}$ by noting that for each $A_\ell$ we can take $t = \frac{k+1}{2k-\ell-1}$.
\end{remark}

\subsection{The weighted inertia bound is far from tight for graph powers}\label{sec:inertiagraphpowers}

With this theorem in hand we can now present, for each fixed $k \ge 3$, our construction of graphs $G_n$ that have $\alpha_k(G_n) = 2$ and $n_{\ge 0}(G_n^k) = \tilde{\Omega}(n^{1/(2k)})$. The graphs $G^k$ will be random Cayley graphs, which gives us control over their second largest eigenvalue. 

\begin{theorem}[\cite{alonRoichman1994}]\label{thm:randomlambda}
    For every abelian group $\Gamma$ of odd order $n$ and for every $\ell = o(n)$, if $S$ is a set of $\ell$ random members of $\Gamma$, then with high probability $\lambda(A(\mathrm{Cay}(\Gamma,S\cup -S))) = O(\sqrt{\ell\log n})$.
\end{theorem}

\begin{remark}
    Compared to the original statement the additional requirements that $\Gamma$ has odd order and $\ell = o(\sqrt{n})$ imply that the multigraph $X(\Gamma,S)$ in \cite{alonRoichman1994} is the simple graph $Cay(\Gamma, S \cup -S)$ with high probability, because $x \neq -x$ and $P(\exists x : x,-x \in S) \le n\left(\frac{\ell}{n}\right)^2 = o(1)$.
\end{remark}

\begin{theorem}\label{thm:kindepseparation}
    For fixed $k \geq 3$ there is a sequence of graphs $G_n$ of order $N$ such that $\alpha_k(G_n) = 2$ and $n_{\geq 0}(G_n^k) \geq \widetilde{\Omega}(N^{1/(2k)})$.
\end{theorem}

\begin{proof}
    For a fixed $k \geq 3$ we pick some prime $q > k$ and then for each $n$ we consider random sets $S \subseteq (\Z/q\Z)^n \setminus 0$ where each non-zero element is included with probability $p(n) =  D_pN^{(1-k)/k}(\log N)^{-(k-1)/k}$ (where $N = q^n$) for small enough constant $D_p$. With high probability the set $S \cup -S$ will have 
    \begin{itemize}
        \item size $\Theta(N^{1/k}(\log N)^{(1-k)/k})$,
        \item no triples $x,y,z$ with $x+y+z = 0$,
        \item $\lambda(A(\mathrm{Cay}(\Gamma,S \cup -S))) = O(N^{1/2k}\log^{1/2k} N)$.
    \end{itemize}
    
    These conditions follow from the following standard inequalities. Let $Y$ be the random variable that counts the number of triples $x,y,z \in S \cup -S$ such that $x + y + z = 0$. Note that $0 \notin S \cup -S$, so in counting tuples for $Y$ none of $x, y, z$ can be inverses of another. As $q > 3$ not all three can be equal. Hence $x,y,z$ are independently included unless two are equal. There are $N-1$ tuples with $x = y$ as $z = -2x$ in that case. The equalities $x = z$ and $y = z$ hold similarly often. The number of triples where $x,y,z$ are all distinct is at most $N^2$, because two values determine the third. The probability that one element is in $S \cup -S$ is $2p-p^2$, hence $\mathbb{E}(Y) \leq 3N(2p-p^2)^2 + N^2(2p-p^2)^3 = o(1)$, so by Markov's inequality $\mathbb{P}(Y \geq 1) \to 0 $. 
    
    By Chernoff's bound we have $\mathbb{P}(||S|- (N-1)p(N)| \geq \frac12(N-1)p(N)) \leq 2e^{-\frac1{12}Np(N)} \to 0$, so $|S \cup -S|$, which is sandwiched between $|S|$ and $2|S|$ has the right size. The last condition follows from Theorem~\ref{thm:randomlambda} with $|S| = \Theta(N^{1/k}(\log N)^{(1-k)/k})$. 
    
    Pick a random $S$ that satisfies all three conditions. Then $|S| = \Theta(N^{1/k}(\log N)^{(1-k)/k})$ and we can choose $D_p$ small enough that Theorem~\ref{thm:alonextendsym0} applies with $F := S \cup -S$. 
        Let $A$ be the symmetric set from the theorem. Choose $G_n = \mathrm{Cay}((\Z/q\Z)^n, A\setminus 0)$. As $0 \in A$ we find that $A \subseteq A+A$ and hence $G_n^k = \mathrm{Cay}((\Z/q\Z)^n, kA\setminus 0)$, which implies that $\bar{G_n^k} = \mathrm{Cay}((\Z/q\Z)^n, F)$. This is a triangle-free graph, because $F$ has no triples summing to zero. It is also a $(N, |F|, O(\sqrt{|F|\log N}))$-graph. Theorem~\ref{thm:expanderinertia} applies, so $G_n$ is a graph with $\alpha_k(G_n) = 2$  and $n_{\geq 0}(G_n^k) \geq \frac{d + \lambda}{2\lambda} = \widetilde{\Omega}\left(|F|^{1/2}\right) = \widetilde{\Omega}(N^{1/(2k)})$.
\end{proof}

\begin{remark}
    For $k = 2$ we cannot use the same $p(N)$, because the triangle-free condition requires $p(N) = o(N^{-2/3})$. The same argument with $q > 3$ and $p(N) = N^{-2/3}/h(N)$, where $h$ is a function that goes to infinity very slowly, gives $n_{\ge 0}(G_n^2) \ge \Omega\left(\left(\frac{N^{1/3}}{h(n)\log N}\right)^{1/2}\right)$ which is only slightly better than our bound for $k = 3$.
\end{remark}

\section{Tightness of the inertia-type bound}\label{sec:inertiatype}

The inertia-type bound is a bound that is computationally simpler to optimize than the weighted inertia bound. However, in return the bound is less strong \cite{notekindep2026}. In this section we prove extreme separation, where the inertia-type bound grows linearly in the order and the $k$-independence number is $1$, i.e. the graph power is complete. We show this separation for large random graphs. One may suspect that the bad performance comes from the non-constant diagonal of $p(A)$, or more generally, the lack of symmetries, in these random graphs. However, we also give an example of vertex-transitive graphs with the same extreme separation. Furthermore, this construction is explicit. Both the random graphs and our explicit graphs have $\alpha_k(G) = 1$ for $k \ge 2$, which even the unweighted inertia bound on $G^k$ can recognise, but the inertia-type bound on $G$ cannot. 

Roughly speaking, our graphs must satisfy two conditions: every vertex is reachable in at most $k$ steps, while the normalized spectral distribution has large enough support to force a linear amount of eigenvalues of both signs for all polynomials. We establish the latter through convergence of the spectral distribution to distributions with infinite support.

In particular for random graphs we mention Wigner's semicircle law: the probability distribution on $[-2,2]$ with density function $\frac1{2\pi}\sqrt{4-x^2}$ and associated distribution function $F_{sc}$ and measure $\mu_{sc}$. For a multiset $B$ the \emph{empirical distribution} is the distribution of a random element. For a Hermitian $n\times n$-matrix $M$ the \emph{Empirical Spectral Distribution} (ESD) is the probability distribution with distribution function $F_M(x) :=\frac1{n}|\{i : \lambda_i \le x\}|$.

Recall that $P_k$ denotes the vector space of real polynomials of degree at most $k$ and that $\In_k(p,G)$ represents the inertia-type bound using a polynomial $p \in P_k$. If $X$ is a random variable which uniformly picks an eigenvalue of $G$, then for all $p \in P_k$,

\begin{align*}
    \In_k(p, G) &= \min\{ |\{i : p(\lambda_i) \le W(p)\}|,\, |\{i : p(\lambda_i) \ge w(p)\}| \} \\ &= n \cdot\min\{ \mathbb{P}(p(X) \le W(p)) , \,\mathbb{P}(p(X) \ge w(p))\}.
\end{align*}

The average diagonal value of $p(A)$ is $\frac1n\tr p(A)$. It is well-known that the trace is equal to the sum of the eigenvalues, so $w(p) \leq \frac1n\sum p(\lambda_i ) = \mathbb{E}[p(X)] \leq W(p)$ and therefore 
\begin{align}\label{eq:inertiatypeprob}
    \In_k(p, G) \ge n \cdot\min\{ \mathbb{P}(p(X) \le \mathbb{E}[p(X)]) , \,\mathbb{P}(p(X) \ge \mathbb{E}[p(X)])\}.
\end{align}

We bound these probabilities for certain sequences of converging measures and show that the ESDs form such a sequence.

\begin{lemma}
\label{lem:positiveprobability_adapted}
Fix an integer $k \ge 1$. Let $\mu$ be a probability measure supported on a compact interval $I \subseteq \mathbb{R}$, with $\mu(S) < 1$ for any set $S \subseteq I$ of size at most $k$. Let $(\mu_n)_{n=1}^\infty$ be probability measures on
$\mathbb R$. Define
\[
    S_k:=\{r\in P_k:\|r\|_{\infty,I}=1\},
    \qquad \text{where }
    \|r\|_{\infty,I}:=\max_{x\in I}|r(x)|.
\]
Assume that
\[
    \sup_{r\in S_k}
    \left|
        \int_{\mathbb R}|r|\,d\mu_n-\int_I |r|\,d\mu
    \right|
    \to 0
\]
and
\[
    \sup_{r\in S_k}
    \int_{\mathbb R\setminus I}|r|\,d\mu_n
    \to 0 .
\]
Then there exists a constant $\varepsilon=\varepsilon(k,\mu, I)>0$ and an
integer $N$ such that, for all $n>N$ and all $p\in P_k$, if
$X_n\sim \mu_n$, then
\[
    \mathbb P\bigl(p(X_n)\ge \mathbb E[p(X_n)]\bigr)\ge \varepsilon \qquad \text{and} \qquad \mathbb P\bigl(p(X_n)\le \mathbb E[p(X_n)]\bigr)\ge \varepsilon .
\]
\end{lemma}

\begin{proof}
The functions $\|r\|_{\infty,I}$ and
\[
    \|r\|_{1,\mu}:=\int_I |r|\,d\mu
\]
are norms on $P_k$, because $\mu(S) < 1$ on the roots of any nonzero polynomial of degree at most $k$.
Since $P_k$ is finite-dimensional and $I$ is compact these norms are equivalent on $P_k$. The equivalence of norms means that there exists $c > 0$ such that  $\|r\|_{1,\mu} \geq c \|r\|_{\infty,I}$ for every $r \in P_k$. The value of this $c$ is exactly
\[
    c:=\inf_{r\in S_k}\int_I |r|\,d\mu >0 .
\]
By the two assumptions, for all sufficiently large $n$,
\begin{equation}\label{eq:convergenceofmeasures}
       \sup_{r\in S_k}
    \left|
        \int_{\mathbb R}|r|\,d\mu_n-\int_I |r|\,d\mu
    \right|
    \le \frac c4
\end{equation}

and
\begin{equation}\label{eq:outside}
       \sup_{r\in S_k}
    \int_{\mathbb R\setminus I}|r|\,d\mu_n
    \le \frac c8 .
\end{equation}
 
We claim that $\varepsilon=c/4$ works.

Fix $p\in P_k$. If $p$ is constant, the claim is immediate. Otherwise, set $m_n:=\int_{\mathbb R}p\,d\mu_n$, which is finite because $p$ is bounded on $I$ and because \eqref{eq:outside} holds, and define the polynomial
\[
    r_n:=\frac{p-m_n}{\|p-m_n\|_{\infty,I}},
\]
which can be viewed as a normalization of $p$.
Then $r_n\in S_k$, so writing $r_n^+=\max(r_n,0)$ we have $0 \leq r_n^+\le 1$ on $I$ and thus
\[\mathbb{P}(p(X_n) \geq m_n) = \mathbb{P}(r_n(X_n) \geq 0)\geq \int_I \mathbf{1}_{\{r_n \geq 0\}} d\mu_n \ge \int_I r_n^+ d\mu_n\]
By definition of $m_n$ we get that
$
    \int_{\mathbb R} r_n\,d\mu_n=0,
$
so
\[
    \int_{\mathbb R} r_n^+\,d\mu_n
    =
    \frac12\int_{\mathbb R}|r_n|\,d\mu_n .
\]
Moreover by \eqref{eq:convergenceofmeasures},
\[
    \int_{\mathbb R}|r_n|\,d\mu_n
    \ge
    \int_I |r_n|\,d\mu-\frac c4
    \ge
    \frac{3c}{4}.
\]
Hence
\[
    \int_{\mathbb R} r_n^+\,d\mu_n
    \ge
    \frac{3c}{8}.
\]
The contribution from outside $I$ is at most the contribution of $|r_n|$ which is at most $c/8$ by \eqref{eq:outside}, so
\[
    \int_I r_n^+\,d\mu_n
    \ge
    \frac{3c}{8}-\frac c8
    =
    \frac c4 .
\]
Since $0\le r_n^+\le 1$ on $I$, it follows that $\mathbb P(r_n(X_n)\ge 0)\ge \frac c4$.
Equivalently,
\[
    \mathbb P\bigl(p(X_n)\ge m_n\bigr)\ge \frac c4.
\]
Applying the same argument to $-p$ gives $\mathbb P\bigl(p(X_n)\le m_n\bigr)\ge \frac c4$.
This proves the lemma.
\end{proof}

\begin{remark}
     When $\mu = \mu_{sc}$ is the semicircle law, then we can show that $\varepsilon \ge \frac{\sqrt{3}}{30\pi k^3}$.
     In \cite{Tikhanovtseva2014} the author studied for $r \in P_k$ with $\int_{[-1,1]} r d\nu = 0$ and $\nu = (1-x^2)^\alpha$ the  value of $$ \min_{r} \int_{\chi(r)} d\nu$$ where $\chi(r) = \{x \in [-1,1] : r(x) \ge 0\}$. Tikhanovtseva showed that this minimum grows as $\Theta(k^{-2(\alpha+1)})$ as a function of $k$. For $\alpha = \frac12$ this gives the minimum for (a scaled) semicircular law. This does not immediately apply to Lemma~\ref{lem:positiveprobability_adapted}, because we have a sequence of measures converging to $\mu$. However, a proof in that paper does suggest a lower bound on $c$. See the proof of Lemma 2 in \cite{Tikhanovtseva2014} for the details. Fix a polynomial $r \in S_k$ and let $t_2$ be a point such that $|r(t_2)| = 1$. Use Markov's inequality on polynomials to lower bound $|r(t)| \ge 1 - k^2|t-t_2|$. Following the same chain of inequalities we get 
     \begin{align*}\int_{[-2,2]} |r| d\mu_{sc} \ge \int_{[-2,2] \cap [t_2-\frac1{k^2}, t_2 + \frac1{k^2}]} (1-k^2|t-t_2|) d\mu_{sc} &\ge \frac{1}{2\pi}\int_{[0,\frac1{k^2} ]} (4x-x^2)^{1/2}(1-k^2x)dx\\ &\ge \frac{\sqrt{3}}{2\pi(1/2+1)(1/2+2)}\frac{1}{k^{2(1/2+1)}}.
     \end{align*}Hence $c \ge \frac{4\sqrt{3}}{30\pi}k^{-3}$.
\end{remark}

\subsection{The optimal inertia-type bound is far from tight for random graphs}\label{sec:inertiatypefar}

 Many models of random Hermitian matrices have spectral distribution converging to the semicircle law in probability, see e.g.\  \cite{intro_randommatrices_book}. This also holds for the normalized adjacency matrix of a $G(n,p)$ graph. This was stated and shown explicitly in \cite{tran_sparse_2013}. Note that the authors do not explicitly mention that the convergence happens `in probability', but that this is the type of convergence present in Wigner's original result. 

\begin{theorem}[{\cite[Theorem~1.3]{tran_sparse_2013} }]\label{thm:semicirclelawgnp}
    Consider the random graphs $G(n,p)$ where $np(1-p) \to \infty$ as $n \to \infty$. 
    Let $A(G)$ denote the adjacency matrix of the graph $G$. The empirical spectral distribution  of the matrix $$\frac{1}{\sqrt{np(1-p)}}A(G(n,p))$$
converges in probability to the semicircle distribution.
\end{theorem} 

Much stronger local semicircle laws have been proven, see e.g.\ \cite{lectureslocallaw}, but the convergence on global scale suffices for our purposes. 

\begin{corollary}
\label{cor:inertiatypefar_sparse}
Fix $k\ge 2$. Let $q_n$ be a sequence of edge probabilities such that as $n \to \infty$, with $d_n:=nq_n$,
\[
    q_n\to0,
    \qquad
    d_n\to\infty,
    \qquad
    \frac{d_n^k}{n} = \omega(\log n),
    \qquad
    d_n=o(n^{2/k}).
\]
Let $G_n=G(n,q_n)$. Then, with high probability,
\[
    \alpha_k(G_n)=1
    \qquad\text{and}\qquad
    \In_k(G_n)=\Omega(n).
\]
In particular, this holds for $q_n=n^{-1+\frac{3}{2k}}$.
\end{corollary}

\begin{proof}
Let $\sigma_n:=\sqrt{nq_n(1-q_n)}$
and $M_n:=\sigma_n^{-1}{A(G_n)}$.
Denote the eigenvalues of $M_n$ by $\xi_1,\ldots,\xi_n$, and let $\mu_n$ be their empirical spectral distribution.

We use two standard spectral facts for $G(n,q_n)$ in this regime. First,
$\mu_n$ converges weakly in probability to the semicircle law $\mu_{sc}$, because $d_n \to \infty$, see Theorem~\ref{thm:semicirclelawgnp}.
Second, applying Lemma~4.3 and Remark~2.4 of \cite{localsemi13} to $M_n - \frac{q_n}{\sigma_n}(J-I)$, says that all non-Perron eigenvalues of $M_n$ lie in
a fixed compact interval, say $[-T,T]$ with $T$ slightly larger than $2$ for large $n$. It also follows that the Perron eigenvalue satisfies $|\xi_1|=O(\sqrt{d_n})$.

Let $I=[-2,2]$. The semicircle measure $\mu_{sc}$ is supported on $I$
and has density strictly positive almost everywhere on $I$.

We verify the hypotheses of Lemma~\ref{lem:positiveprobability_adapted}.
Let $S_k:=\{r\in P_k:\|r\|_{\infty,I}=1\}$.
Since $S_k$ is compact in the finite-dimensional space $P_k$, there is a
constant $C_k>0$ such that, for all $r\in S_k$,
\[
    |r(x)|\le C_k(1+|x|^k)
    \qquad\text{for all }x\in\mathbb R.
\]
Therefore the Perron eigenvalue contributes at most
\begin{equation}\label{eq:perroncontribution}
    \frac1n \sup_{r\in S_k}|r(\xi_1)|
    \le
    O\left(\frac{1+d_n^{k/2}}{n}\right)
    =
    o(1)
\end{equation}
to $\int_{\mathbb{R}} |r| d\mu_n$, because $d_n=o(n^{2/k})$.

On the remaining $n-1$ eigenvalues, the functions $|r|$, $r\in S_k$,
are uniformly bounded and form a compact family of continuous functions on
$[-T,T]$. Hence weak convergence in probability to $\mu_{sc}$ implies,  uniformly in
$r\in S_k$, for all fixed $\delta > 0$
\[
    \mathbb{P}\left(\Big|\int_{\mathbb{R}}|r|\,d\mu_n - \int_I |r|\,d\mu_{sc}\Big| \ge \delta \right) \to 0.
\]
This is the right convergence for our argument. Fixing $\delta = c/4$, where $c$ comes from the proof of Lemma~\ref{lem:positiveprobability_adapted} we see that \eqref{eq:convergenceofmeasures} holds with high probability. The measure $\mu_{sc}$ is supported on $I=[-2,2]$, has no atoms and $|r|$ is bounded outside $I$ except at $\xi_1$, so combining \eqref{eq:perroncontribution} and weak convergence, the total
$|r|$-weighted mass outside $I$ is also negligible (with high probability) uniformly in $r\in S_k$:
\[
    \sup_{r\in S_k}
    \int_{\mathbb R\setminus I}|r|\,d\mu_n
    \to 0.
\]
Thus Lemma~\ref{lem:positiveprobability_adapted} applies.

Consequently, there exists $\varepsilon>0$, depending only on $k$, such
that with high probability, for every $r\in P_k$,
\[
    \mathbb P\bigl(r(X_n)\ge \mathbb E[r(X_n)]\bigr)\ge \varepsilon
\]
and
\[
    \mathbb P\bigl(r(X_n)\le \mathbb E[r(X_n)]\bigr)\ge \varepsilon,
\]
where $X_n$ is a uniformly chosen eigenvalue of $M_n$.

Now return to the unnormalized adjacency eigenvalues. If
$\lambda_i=\sigma_n\xi_i$, then for every $p\in P_k$ the polynomial
\[
    r(x):=p(\sigma_n x)
\]
also belongs to $P_k$. Hence the preceding probability bounds apply to
$p(\lambda_i)$, uniformly over all $p\in P_k$. Therefore by \eqref{eq:inertiatypeprob}
\[
\begin{aligned}
    \In_k(p,G_n)
    &\ge
    n\cdot
    \min\left\{
        \mathbb P\bigl(p(X)\le \mathbb E[p(X)]\bigr),
        \mathbb P\bigl(p(X)\ge \mathbb E[p(X)]\bigr)
    \right\} 
    \ge
    \varepsilon n,
\end{aligned}
\]
where $X$ is a uniformly chosen eigenvalue of $A(G_n)$. Since this holds
for all $p\in P_k$, we get $\In_k(G_n)=\Omega(n)$
with high probability.

It remains to show that $\alpha_k(G_n)=1$. By the standard diameter
threshold for binomial random graphs, if
\[
    \frac{d_n^k}{n} =  \omega(\log n),
\]
then $G(n,q_n)$ has diameter at most $k$ with high probability, see \cite[Chapter~10]{Bollobas_2001}. This is
one of our assumptions. Hence, with high probability, every two vertices of
$G_n$ are at distance at most $k$, and therefore $\alpha_k(G_n)=1$.
The proof is complete.

Finally, for the explicit choice $q_n=n^{-1+3/(2k)}$, we have $d_n=nq_n=n^{3/(2k)}$, so
\[
    \frac{d_n^k}{n}=n^{1/2} = \omega( \log n)
\]
and
\[
    \frac{d_n^{k/2}}{n}=n^{-1/4}\to0.
\]
Thus all assumptions above are satisfied.
\end{proof}

\subsection{An explicit vertex-transitive graph class with far from tight inertia-type bound}\label{sec:cayleyinertiatypefar}

A graph is \emph{vertex-transitive} when for any pair of vertices $u,v$ there is an automorphism sending $u$ to $v$. In this section we construct an explicit vertex-transitive graph that obtains the same separation as the random graphs in Section~\ref{sec:inertiatypefar}. When $G$ is a vertex-transitive graph, such as a Cayley graph, the diagonal elements $(p(A))_{ii}$ in the inertia-type bound from \Cref{thm:inertia-type} are all equal and hence rephrasing using the empirical spectral distribution $X$ and $Y = p(X)$, the inertia-type bound states
\[\alpha_k(G) \leq n\cdot \min\{\mathbb{P}(Y \geq \mathbb{E}[{Y}]),\mathbb{P}(Y \leq \mathbb{E}[{Y}])\}. \]

In this section we prove the following result. 

\begin{theorem}\label{thm:cayleyinertiatypefar}
    Fix an integer $ k \ge 5$ and prime $p > k$. Define $s = \lfloor \frac{k}2\rfloor +1$. There exists an integer $a \ge 1$ such that $s \mid p^{2a} -1$. For any large enough prime power $q = (p^{2a})^h$ the graph 
    \[
    G_q^{s} := \mathrm{Cay}(\F_q^{s}, S_{s}), \quad \text{ where } S_{s} = \{\pm(t, t^2, \ldots, t^s) : t \in \F_q^\times\}
    \]
    has
    \[\alpha_k(G_q^{s}) = 1, \quad \quad \text{and} \quad \quad \In_k(G_q^{s}) = \Omega(|V(G_q^{s})|) \text{ as } q \to \infty.\]
    
\end{theorem}

\begin{remark}
    We should note that a similar result holds for $k = 2,3,4$, but the proofs given in this section work assuming $k \ge 5$. We outline here shortly what needs to change for the small $k$. For $k = 4$ we can use the same families of graphs, but we need a different proof that the diameter is still at most $k$. For $k = 2,3$ we have to rely on different graphs related to Kloosterman and Birch sums which also satisfy equidistribution. For $k = 3$ we can use $S = \{(t,t^3) : t \in \F_q^\times\}$ and rely on equidistribution of Birch sums as stated in \cite{FFKW25}. It can be shown that these graphs have diameter $3$. For $k = 2$ one can use the `Euclidean' Cayley graph with connection set $\{(x,y,z,w) \in \F_q^4 : x^2 + y^2+z^2+w^2 = 1\}$, see \cite{brittenham25} for the spectral properties and \cite{lorenz25} for the diameter.  
\end{remark}

For a $d$-regular graph the \textit{normalized spectrum} $\mathrm{NSpec}(G)$ is obtained from $\mathrm{Spec}(G)$ by dividing each element by $\sqrt{d}$. 

Inspired by the work in \cite{FFKW25} we seek a family of graphs defined by algebraic relations that have a spectrum converging to a distribution with infinite support.

In particular, we let $p$ be an odd prime and $q = p^h$ a power of $p$. We consider the Cayley graph $G_q^s = \mathrm{Cay}(\F_q^s, S)$ for $s \geq 3$ where $S = \{\pm(t, t^2, \ldots, t^s) : t \in \F_q^\times\}$. This graph has order $n = q^s$ and degree $2(q-1)$. We show that the non-principal eigenvalues are of order at most $\sqrt{q}$ and that the normalized spectrum converges to a distribution with infinite support. The same argument as found in Corollary~\ref{cor:inertiatypefar_sparse} then implies that $\In_k(G_q^s) \ge \Omega(n)$. Additionally, we show that for $s = \lfloor \frac{k}{2}\rfloor + 1$ this graph has diameter at most $k$ which implies that $\alpha_k(G_q^s) = 1$.

For the diameter we use Cohen's result that all factorization patterns of polynomials are almost equally common \cite{cohen1972}. For a partition $\lambda$ of $n$ let $T(\lambda)$ be the proportion of elements in the symmetric group with cycle pattern $\lambda$. For a field $\F$ and a polynomial $P \in \F[x]$ the \emph{factorization pattern} is the multiset of degrees of the irreducible factors of $P$.

\begin{theorem}[\cite{cohen1972}]\label{thm:cohen}
Let $A$ be a set consisting of all polynomials in $\F_q[x]$ of degree $d$ with some set of $s$ coefficients preassigned (the constant term being non-zero, if it is preassigned). Suppose that $A$ is not a subset of $\F_q[x^\ell]$ for any $\ell > 1$. If $q = p^h$ for some prime $p > d$ and $s \le d - 1$, then, for any factor pattern $\lambda$, the number of polynomials with factor pattern $\lambda$ is
$$T(\lambda)q^{d-s+1} + O(q^{d - s+\frac12}),$$
where the implied constant depends only on $d$. 
\end{theorem}

In the case that $\lambda = 1^{(d)}$ this theorem counts the number of polynomials that split fully into linear terms.

\begin{lemma}\label{lem:curve_diameter}
    For fixed $s \le k - 2$, fixed prime $p > k$ and large enough $q = p^h$, all $x \in \F_q^s$ can be expressed as a sum of $k$ elements from the set $S_+ = \{(t, t^2, \ldots, t^s) : t \in \F_q\}$.
\end{lemma}

\begin{proof} 
    Consider an arbitrary element $x = (x_1,\ldots, x_s) \in \F_q^s$. 
    We look for  $t_1,\ldots, t_k \in \F_q$ such that the equations
    \begin{equation}\label{eq:diameter} 
    \sum_{i = 1}^k t_i^j = x_j, \qquad \text{for all }1 \leq j \leq s
    \end{equation}
    hold. Let $e_j$ denote the $j$th elementary symmetric function of the $t_i$, i.e., $e_1 = t_1+\ldots +t_k$, $e_2 = \sum_{j_1 < j_2} t_{j_1}t_{j_2}$ etc. The power sums and the elementary symmetric functions are related by the Newton identities which are invertible because $p > k$, so the equations in \eqref{eq:diameter} determine the values of $e_1,\ldots, e_s$. Let $A$ be the set of monic degree $k$ polynomials in $\F_q[x]$ with the first coefficients being the fixed values $1, -e_1, e_2,\ldots, (-1)^se_s$. This is a set of polynomials with $s+1 \leq k-1$ coefficients fixed and in particular with the coefficient of $x$ not fixed. Hence $A\not\subseteq \F_q[x^\ell]$ for any $\ell > 1$ and so by Theorem~\ref{thm:cohen} the number of polynomials in $A$ with factor pattern $1^{(k)}$ is $T(1^{(k)})q^{k-s} + O(q^{k-s-\frac12})$. For large enough $q$ this is positive. Hence there is a polynomial $r$ in $A$ that splits into linear terms. The values of $e_1,\ldots, e_s$ determine $\sum t_i^j$ for $j \leq s$. By definition of $A$ the roots of $r$ solve the equations in \eqref{eq:diameter}, so there are $k$ elements of $S_+$ that sum to $x$. 
\end{proof}

For the limiting spectrum we rely on Weil's bound of character sums and Katz' equidistribution theorem. In particular, we use the following formulations of these results.

\begin{theorem}[Weil bound \cite{lidl_finite_1996}] \label{thm:weilbound}
    Let $f \in \F_q[x]$ be a polynomial of degree $d \geq 1$ with $\gcd(d,q) = 1$ and let $\chi$ be a nontrivial additive character of $\F_q$. Then
    \[
    \left|\sum_{a \in \F_q} \chi(f(a))\right| \le (d-1)\sqrt{q}
    \]
\end{theorem}

\begin{theorem}[Our version of Katz' equidistribution theorem]\label{thm:ourequidistribution}
    Let $d \ge 3$ be a fixed natural number and let $ p > 5$ be a prime not dividing $d$. Choose an integer $a \ge 1$ such that $d \mid p^{2a} - 1$ and put $Q = p^{2a}$. Let $q = Q^r$, $\chi$ a nontrivial additive character of $\F_q$ and let $\mathcal{P}_d^0(q) \subseteq \F_q[x]$ be the class of polynomials of degree $d$ with zero constant term. As $r \to \infty$ the empirical distribution associated with the multiset $$\Big\{q^{-1/2} \mathfrak{Re} \left(\sum_{t \in \F_q} \chi(f(t)) \right): f \in \mathcal{P}_d^0(q) \Big\}$$ converges to a measure on $[1-d, d-1]$ with infinite support.
\end{theorem}

This statement follows from deep algebraic geometric results by Katz \cite{KatzSarnak1999, Katz2005}. The derivation of this precise form is done in Appendix~\ref{sec:appendix}. Using these results we prove Theorem~\ref{thm:cayleyinertiatypefar}.

\begin{proof}[Proof of Theorem~\ref{thm:cayleyinertiatypefar}]
For $k \ge 5$ and $s = \lfloor \frac{k}{2} \rfloor + 1$ we have $3 \le s \le k - 2$. Pick a prime $p > k$. Then $p$ is a unit modulo $s$, so there is an $a$ such that $s \mid p^{2a} - 1$. Consider the graph $G_q^s$ for powers $q$ of $p^{2a}$. Lemma~\ref{lem:curve_diameter} states that any element of $\F_q^s$ can be written as the sum of at most $k$ (leaving out the $t = 0$ elements) elements of $S$, so $G_q^s$ has diameter $\leq k$ for large enough $q$ and hence $\alpha_k(G_q^s) = 1$. We show that the conditions of Lemma~\ref{lem:positiveprobability_adapted} are satisfied. For any fixed nontrivial character $\chi$ of $\F_q$ the eigenvalues of the Cayley graph are
\[
\lambda_a = \sum_{u \in S} \chi(a \cdot u).
\]
ranging over all $a \in \F_q^s$. As $S$ consists of $\pm(t,\ldots, t^s)$ and $\chi(-x) = \bar{\chi(x)}$ this sum becomes $2\sum_{t \in \F_q^\times} \mathfrak{Re}\,  \chi(a_1t+a_2t^2+\ldots+a_st^s)$. Theorem~\ref{thm:weilbound} applies, because $p > s \ge \deg f$, so $\gcd(p, \deg f) = 1$. Additionally, using $\chi(0) = 1$ the eigenvalues for $a \neq 0$ fall in $[-2(s-1)\sqrt{q} - 2, 2(s-1)\sqrt{q} - 2]$. The Perron eigenvalue, the degree, is $2(q-1)$. Let $\mu_q$ denote the ESD of the normalized adjacency matrix $\frac{1}{\sqrt{2(q-1)}}A(G_q^s)$ whose eigenvalues are denoted $\xi_a$. All the nontrivial normalized eigenvalues fall in the interval $I = [-\sqrt{2}(1+s),\sqrt{2}(s-1)]$. For the inertia-type bound we are interested in $r(\xi_a)$ for all $r \in S_k = \{f \in P_k : ||f||_{\infty,I} = 1\}$. We have $\xi_0 = (1+o(1))\sqrt{2q}$ and because the order $n$ is $q^s$, we have
\[
\frac1n\sup_{r \in S_k} |r(\xi_0)| = O\left(\frac{1+(2q)^{k/2}}{n}\right) = o(1).
\]
 We conclude that $\sup_{r \in S_k}\int_\mathbb{R\setminus I} |r| d\mu_q \to 0$. For the $a \in \F_q^s$ with $a_s \neq 0$ the polynomial $p_a = a_1t +\ldots +a_st^s$ is in $\mathcal{P}_s^0(q)$. The associated eigenvalue $\xi_a$ is equal to
 $$\frac{2}{\sqrt{2(q-1)}}\sum_{t \in \F_q^\times} \mathfrak{Re}\,  \chi(p_a(t)).$$
As $q \to \infty$ the distribution of these eigenvalues differs by $o(1)$ from 
 $$\Big\{\sqrt{2}q^{-1/2} \mathfrak{Re} \left(\sum_{t \in \F_q} \chi(f(t)) \right): f \in \mathcal{P}_s^0(q) \Big\},$$
 which converges to some measure $\mu$ with infinite support (contained in $I$) by Theorem~\ref{thm:ourequidistribution}. The contribution of the eigenvalues with $a_s = 0$ goes to zero. This was already shown for the Perron eigenvalue. The remaining eigenvalues lie in $I$, so $|r(\xi)| \le 1$ for $r \in S_k$. The eigenvalues with $a_s = 0$ form a $\frac1q$ fraction of all eigenvalues, so their contribution to $\int_\mathbb{R} |r| d\mu_q$ goes to zero. Therefore the whole ESD converges to $\mu$ uniformly over $r$, i.e., 
\[
\sup_{r\in S_k} \left|\int_\mathbb{R} |r| d\mu_q - \int_I |r| d\mu\right| \to 0.
\]
This means that all conditions of Lemma~\ref{lem:positiveprobability_adapted} are satisfied and hence there is an $\varepsilon > 0$, such that the inertia-type bound is at least $\varepsilon q^s$. 
\end{proof}

\section{Concluding remarks}

It is now clear that the weighted inertia bound, and even more so the inertia-type bound, can be far from tight and that this is not an exception, but happens for almost all graphs of the right density. Still, there are examples where these bounds are effective in Erd\H{o}s-Ko-Rado results \cite{Calderbank_Frankl_1992} and coding theory \cite{APR2025}. Understanding what properties of the graph or spectrum make the inertia bound perform ``well'' remains an open problem. In this process one may also ask if current separation results are really best possible.

For $k$-independence specifically, we do not claim that our separation results are optimal. An obvious example where improvement should be possible is in Theorem~\ref{thm:kindepseparation}, which does not give a significantly better lower bound for $k = 2$ compared to $k = 3$. Note further that our results give separation for dense  graph powers. For sparse graphs, such as sparse random graphs or random regular graphs the $k$-independence number is larger, see e.g.\ \cite{atkinsonfrieze2004, beislarge2005}, but it also seems likely that the weighted inertia bound performs worse. Currently, it is unknown how much worse the weighted inertia bound (and also the inertia-type bound) will be for sparse graphs if we cannot guarantee $C_4$-freeness. 

There are more lower bounds on the inertia(-type) bound. It is known that the quantum independence number $\alpha_q(G)$, and even the quantum $k$-independence number $\alpha_{kq}(G) := \alpha_q(G^k)$, is in between the independence number and the inertia bounds, for the definition and this result see \cite{WEA2022}. Any graph with separation between $\alpha(G)$ and $\alpha_q(G)$ also gives separation between $\alpha(G)$ and $n_{\ge 0}(G)$. However, it is difficult to determine $\alpha_q$ and hence to give large separation between $\alpha(G)$ and $\alpha_{q}(G)$ or $\alpha_k(G)$ \ and $\alpha_{kq}(G)$. In particular, it is known that $\alpha_q = \alpha$ if $\alpha(G) \le 2$. Therefore the results in Theorems~\ref{thm:kindepseparation} and \ref{thm:cayleyinertiatypefar} and Corollary~\ref{cor:inertiatypefar_sparse} also give separation between the inertia(-type) bound and $\alpha_{kq}(G)$. Our results do not help in identifying which graphs have (large) separation between the ($k$-)independence number and the corresponding quantum parameter, but we know that those graphs must also have an inertia(-type) bound that is far from tight.

\subsection*{AI tools declaration}

The research of this paper was assisted by the use of Gemini and ChatGPT. Specifically, the ideas for considering random graphs, random Cayley graphs and pseudorandom Cayley graphs were conceived by the authors, but specific proofs and constructions for the latter two were suggested by AI. The working out of details of the proof in Section~\ref{sec:inertiagraphpowers} was assisted by ChatGPT 5.5. More notably, the proof method of Lemma~\ref{lem:many_sets} was suggested by Gemini 3.1. Additionally, Section~\ref{sec:cayleyinertiatypefar} and Appendix~\ref{sec:appendix} were greatly supported by ChatGPT 5.5 and 5.6. The construction used in Theorem~\ref{thm:cayleyinertiatypefar} was suggested by ChatGPT 5.5 based on our specific suggestions for such a construction. ChatGPT 5.5 and 5.6 were consequently used for understanding the algebraic geometric background of these constructions. These tools were also used throughout to find literature, typos and small mistakes. However, all writing in this paper was done by the authors.

\subsection*{Acknowledgements}

Aida Abiad and Nils van de Berg are supported by NWO (Dutch Research Council) through the grant  VI.Vidi.213.085  (grant ID \url{https://doi.org/10.61686/CSWNV56135}). Sam Mattheus is supported by postdoctoral fellowships 1267923N and 12A0J27N from the Research Foundation Flanders (FWO). This research was partially supported by the FWO Scientific Research Communities: Graphs, Association schemes and Geometries: structures, algorithms and computation (W003324N) and Finite Geometry, Coding Theory and Cryptography (W0.032.24N). The authors thank Harper Reijnders for useful discussions on Section \ref{sec:weightedinertia} and Prof. Emmanuel Kowalski for giving feedback on the arithmetic geometry content of the appendix. Part of this work was done while the first author was visiting the Simons Institute for the Theory of Computing at UC Berkeley.

\printbibliography[heading=bibintoc]


\newpage

\appendix

\section{Appendix: Proof of Theorem~\ref{thm:ourequidistribution}}\label{sec:appendix}

Recall that we want to prove the following.

 \newtheorem*{thm:appendix}{Theorem \ref{thm:ourequidistribution}}
\begin{thm:appendix}[Our version of Katz' equidistribution theorem]
    Let $d \ge 3$ be a fixed natural number and let $ p > 5$ be a prime not dividing $d$. Choose an integer $a \ge 1$ such that $d \mid p^{2a} - 1$ and put $Q = p^{2a}$. Let $q = Q^r$, $\chi$ a nontrivial additive character of $\F_q$ and let $\mathcal{P}_d^0(q) \subseteq \F_q[x]$ be the class of polynomials of degree $d$ with zero constant term. As $r \to \infty$ the empirical distribution associated with the multiset $$\Big\{q^{-1/2} \mathfrak{Re} \left(\sum_{t \in \F_q} \chi(f(t)) \right): f \in \mathcal{P}_d^0(q) \Big\}$$ converges to a measure on $[1-d, d-1]$ with infinite support.
\end{thm:appendix}

This result nearly follows from \cite{entin}, specifically Section~2 up to Theorem~3 where Entin states the following as a special case of Katz' results.

\begin{theorem}[\protect{\cite[Theorem~3]{entin}}]
    Let $d \ge 3$ be a fixed natural number and assume $\F_q$ is a field whose characteristic is a fixed prime $> 5$ that does not divide $d$. Then as $q  \to \infty$ the family $\{\Theta_f\}_{f \in \mathcal{P}_d}$ becomes equidistributed in the space of conjugacy classes of some Lie
group $SU(d - 1) \subset G \subset U(d -1)$ with the measure induced from the Haar
measure on $G$.
\end{theorem}

However, this statement is not precise. The group $G$ may change as $q$ varies. This can be resolved  by passing to field extensions of $\F_Q$. We show this essential step in Lemma~\ref{lem:equalitymonodromy}. Pasting this result into Entin's statement does not give a satisfying rigorous proof. Instead, we provide a more complete argument. We have aimed to write this proof so it can be understood with as little algebraic geometric background as possible. However, the argument uses some deep results in arithmetic geometry. We have tried to translate these to our (simpler) setting. For a complete understanding the reader should be familiar with representation theory, Haar measures, the topology of linear groups and the fundamental groups of schemes.

We first state the two deep results as written by Katz and Sarnak that we will be using.

\begin{theorem}[\protect{\cite[Theorem~3.9.2]{Katz2005} for $n = 1$, $p > 5$}]\label{thm:katzmonodromy}
    Let $k$ be a finite field with characteristic $p > 5$ and $\ell \neq p$ both prime, let $\psi : k \to \bar{\mathbb{Q}}_\ell$ be a nontrivial additive character and let $d$ be an integer such that $\gcd(p,d) = 1$. Let $E$ be a finite field extension of $k$ and let $D(d, E) \subseteq E[x]$ denote the scheme of polynomials of degree $d$ over $E$. Then the lisse, pure of weight zero sheaf $M(d,\psi)$ of rank $d-1$ on $D(d,E)$ with trace function given by
    \[
    \tr \rho(Frob_{E,f} \mid M(d,\psi))= -\frac{1}{|E|^{1/2}}\sum_{t \in E} \psi(f(t))
    \]
    for $f \in D(d,E)$ has geometric monodromy group $G_{geom}(\bar{\mathbb{Q}}_\ell)$ contains $SL(d-1, \bar{\mathbb{Q}}_\ell)$.
\end{theorem}

    For a lisse sheaf on scheme $X$, the geometric monodromy group is the Zariski closure of $\rho(\pi_1(X/\bar{k}))$ where $\rho$ is the representation associated to the sheaf, so in our case the group lives in $GL(d-1, \bar{\mathbb{Q}}_\ell)$. We note that the theorem statement is also true when restricting the sheaf to the set of polynomials with zero constant coefficient. This can be seen as \cite[Theorem~3]{entin} applies to this same set of polynomials with zero coefficient. 

\begin{theorem}[\protect{\cite[Theorem~9.2.6]{KatzSarnak1999} for $n = 1$, $p > 5$}]\label{thm:katzequidistribution}
Let $k$ be a finite field, $X/k$ a smooth, geometrically connected scheme, $\ell$ a prime number invertible in $k$, $\iota$ a field embedding of $\bar{\Q}_\ell$ into $\mathbb{C}$
, and $\mathcal{F}$ a lisse
$\bar{Q}_\ell$-sheaf on $X$ of rank $r \ge 1$ corresponding to an $r$-dimensional continuous $\bar{Q}_\ell$-
representation $\rho$ of $\pi_1(X/k)$. Assume
that $\mathcal{F}$ is $\iota$-pure of weight zero, and that $\rho(\pi_1) \subseteq G_{geom}(\bar{\mathbb{Q}}_\ell)$. Denote by $K$ a
maximal compact subgroup of $G_{geom}(\mathbb{C})$. Then the Frobenius conjugacy classes $\theta(E, x)$ are equidistributed in the space of conjugacy classes of $K$ as $|E| \to \infty$, in the sense that for any continuous $\mathbb{C}$-valued function $g$ that is constant on conjugacy classes $ \frac{1}{|X(E)|} \sum_{x \in X(E)} g(\theta(E,x)) \to \int_K g \, d\mu_{Haar}$.
\end{theorem}

In order to use the equidistribution theorem we need to show the following.

\begin{lemma}\label{lem:equalitymonodromy}
    The condition $\rho(\pi_1) \subseteq G_{geom}$ holds in the setting of Theorem~\ref{thm:katzequidistribution} when $k = \F_Q$, $X = \mathcal{P}_d^0$ and $\mathcal{F} = M^0(d,\psi)$, the sheaf from Theorem~\ref{thm:katzmonodromy}.
\end{lemma} 

\begin{proof}
    
    Let $f_a(t) = bt^d$ for $b \in \F_Q^\times$ be a map in $X(\F_Q)$. We will show that there is a choice of $b$ such that $B_b = \rho(Frob_{\F_Q,f_b})$ (a $d-1$ dimensional matrix) has $\det B_b = 1$. Given this, we know that $B_b \in SL(d-1) \subseteq G_{geom}$ by Theorem~\ref{thm:katzmonodromy}. As $X$ is connected, there is a short exact sequence
    \[1 \longrightarrow \pi_1(X/\bar{k}) \longrightarrow \pi_1(X/k) \underset{\tau}{\longrightarrow} \mathrm{Gal}(\bar{k}/k) \longrightarrow 1,
    \]
    see e.g. \cite{stacks-project}.  In our case $k = \F_Q$ and it is known that the Galois group is topologically generated by the image of a Frobenius map \cite{LenstraGalois}. By exactness and compactness this means that $\pi_1(X/\F_Q)$ is topologically generated by $\pi_1(X/\bar{\F}_Q)$ and $Frob_{\F_Q,f_b}$. Applying $\rho$ and noting that $B_b$ is in the topologically closed $G_{geom}$ then implies that $\rho(\pi_1(X/\F_Q)) \subseteq G_{geom}$.
    
    All that is left is to find an $b$ with $\det B_b = 1$. We identify the eigenvalues of $B_b$ by considering $\tr B_b^r$. Note that $Frob_{\F_Q,f_b}^r = Frob_{\F_{Q^r},f_b}$, so that $B_b^r = \rho(Frob_{\F_{Q^r},f_b})$. By Theorem~\ref{thm:katzmonodromy} this has trace 
    $$-\frac{1}{Q^{r/2}}\sum_{t \in \F_{Q^r}} \psi_r(bt^d)$$
    where $\psi_r = \psi \circ \tr_{\F_{Q^r}/\F_Q}$ is the lift of $\psi$ to $\F_{Q^r}$
    Due to $d \mid Q-1$ we know that $\F_Q$ contains all $d$ roots of unity. Therefore the calculation of $-\sum_x \psi_r(bx^e)$  in the proof of \cite[Theorem~6.8.34]{Katz2005} applies and the trace can be rewritten using a nontrivial multiplicative character $\eta$ on $\F_Q$. 

    \begin{equation}\label{eq:trace_evals}
        - Q^{-r/2} \sum_{\eta \neq 1, \eta^d = 1} \bar{\eta}_r(b)\sum_y \eta_r(y)\psi_r(y).
    \end{equation}

    By definition we have $$\eta_r(b) = \eta \circ N_{\F_{Q^r}/\F_Q}(b) = \eta(b^r) = \eta(b)^r,$$ because $b \in \F_Q$. In addition, the Hasse-Davenport lifting relation \cite{Hasse35} says that $$- \sum_{y\in \F_{Q^r}} \eta_r(y)\psi_r(y) = (-1)^r \left(\sum_{y \in \F_Q} \eta(y)\psi(y)\right)^r$$. 

    Plugging this into \eqref{eq:trace_evals} gives
    \begin{equation}\label{eq:traceFrob}
        \tr B^r = \sum_{\eta \neq 1, \eta^d = 1}\left(-Q^{-1/2}\bar{\eta}(b)\sum_y \eta(y)\psi(y)\right)^r.
    \end{equation}

    This identifies the eigenvalues of $B$ as $\lambda_\eta = -Q^{-1/2}\bar{\eta}(b)\sum_y \eta(y)\psi(y)$. Using $\eta(b)\eta^{-1}(b) = 1$ and the identity of Gauss sums \cite{Ireland1982} we find that $\lambda_\eta\lambda_{\eta^{-1}} = \eta(-1)$. For odd $d$ we have $\eta(-1) = 1$ and we can pair up all eigenvalues, so the determinant is
    $$
    \prod_{\eta \neq 1, \eta^d = 1} \lambda_\eta =  1.
    $$
    For even $d$ we get one eigenvalue that is not paired. This corresponds to the nontrivial character $\eta_2$ with $\eta_2^2 = 1$. By the Gauss sum identity we have $\lambda_{\eta_2}^2 = \eta_2(-1)$. As $Q$ is an odd square we find that $-1$ is a square in $\F_Q$, so $\eta_2(-1) = 1$. Hence $\lambda_{\eta_2} \in \{-1,1\}$. This value depends on $b$. In particular, we can pair up the inverse pairs to see there is an $\alpha \in \{-1,1\}$ independent of $b$ such that $\det B_b = \alpha\eta_2(b)$.  By taking $b$ to be a square or non-square we can flip the sign of the determinant. Hence, we can choose a $b$ such that $\det \rho(Frob_{F_Q, f_b}) = 1$. 
\end{proof}

Now we can show Theorem~\ref{thm:ourequidistribution}. We take $k = \F_Q$ (so that extensions $E = \F_q$), $X(\F_q) = \mathcal{P}^0_d(q)$ and $\mathcal{F} = M^0(d, \chi)$. Due to Lemma~\ref{lem:equalitymonodromy} and our comments after Theorem~\ref{thm:katzmonodromy} all the necessary conditions of Theorem~\ref{thm:katzequidistribution} are satisfied. The Frobenius conjugacy class $\theta(E,f)$, defined in 9.2.4 in \cite{KatzSarnak1999}, is exactly the conjugacy class of $K$ which is contained in the (semisimple part of the) conjugacy class of $\iota\rho(Frob_{E,f}\mid \mathcal{F})$. Importantly for us, this implies that $\tr \theta(E,f) = 
- q^{-1/2}\sum_{t \in \F_q} \chi(f(t))$. For any continuous function $g$, the composition $g \circ (-\mathfrak{Re}\tr)$ is continuous and constant on conjugacy classes and thus Theorem~\ref{thm:katzequidistribution} states that there is a compact subgroup $K$ of $G_{geom}(\mathbb{C})$ such that
\[
\frac{1}{q^{d-1}{(q-1)}} \sum_{f \in \mathcal{P}_d^0(q)} g \circ \mathfrak{Re}\left(\frac{1}{q^{1/2}} \sum_{t \in \F_q} \chi(f(t))\right) \to \int_K g \circ (-\mathfrak{Re} \tr) \, d\nu_{Haar}.
\]

From Theorem~\ref{thm:katzmonodromy} we get that $SL(d-1, \bar{\mathbb{Q}}_\ell) \subseteq G_{geom}(\bar{\mathbb{Q}}_\ell) \subseteq GL(d-1, \bar{\mathbb{Q}}_\ell)$. The embedding $\iota$, then gives $SL(d-1,\mathbb{C}) \subseteq G_{geom}(\mathbb{C}) \subseteq GL(d-1,\mathbb{C})$. Maximal compact subgroups of these linear groups are conjugations of $SU(d-1)$ and $U(d-1)$, so $K$ can be chosen to lie between these unitary groups. The function $\Phi : K \to \mathbb{R}$ that maps $M$ to $-\mathfrak{Re}\tr M$ is continuous. Then $\mu = \Phi_*\nu_{Haar}$ is the probability measure that assigns to open $I \subseteq \mathbb{R}$ the value $\mu(I) = \nu(\Phi^{-1}(I))$. As $K\subseteq U(d-1)$, this measure is supported on
$[1-d,d-1]$. For any $y \in \mathrm{Im} \Phi$, take an open neighbourhood $J_y$ of $y$. By continuity $\Phi^{-1}(J_y)$ is a non-empty open subset of $K$ and hence the Haar measure, defined on a compact group, assigns positive measure to it. Hence $\mu(J_y) = \nu(\Phi^{-1}(J_y)) > 0$. This holds for all $y$, so $\mu$ has support $\mathrm{Im} \Phi$.

For every $\theta\in\mathbb R$, the matrix
    $U_\theta=
    \operatorname{diag}(e^{i\theta},e^{-i\theta},1,\ldots,1)$
belongs to $SU(d-1)\subseteq K$, and
    $\Phi(U_\theta)=-(d-3+2\cos\theta)$.
Thus $\operatorname{supp}\mu$ contains the nondegenerate interval
$[1-d,5-d]$, and in particular is infinite. We conclude that the empirical distribution of 
$$\Big\{q^{-1/2} \mathfrak{Re} \left(\sum_{t \in \F_q} \chi(f(t)) \right): f \in \mathcal{P}_d^0(q) \Big\}$$ 
converges to a measure $\mu$ with infinite support.

\end{document}